\documentclass{amsart}
\usepackage{amssymb}
\usepackage{amsfonts}
\usepackage{amscd}
\usepackage{amsmath}
\newtheorem{theorem}{Theorem}[section]
\newtheorem{proposition}[theorem]{Proposition}

\newtheorem{corollary}[theorem]{Corollary}

\theoremstyle{definition}
\newtheorem{definition}[theorem]{Definition}

\newtheorem{remark}[theorem]{Remark}

\newcommand{\refE}[1]{(\ref{E:#1})}
\newcommand{\refS}[1]{Section~\ref{S:#1}}

\newcommand{\refP}[1]{Proposition~\ref{P:#1}}
\newcommand{\refD}[1]{Definition~\ref{D:#1}}

\newcommand{\C}{\ensuremath{\mathbb{C}}}
\newcommand{\N}{\ensuremath{\mathbb{N}}}

\renewcommand{\P}{\ensuremath{\mathbb{P}}}
\newcommand{\Z}{\ensuremath{\mathbb{Z}}}
\newcommand{\gb}{\overline{\mathfrak{g}}}
\newcommand{\laxg}{\overline{\mathfrak{g}}}
\newcommand{\g}{\mathfrak{g}}
\newcommand{\gh}{\widehat{\mathfrak{g}}}
\newcommand{\ga}{\gamma}
\renewcommand{\a}{\alpha}
\renewcommand{\b}{\beta}
\newcommand{\ka}{\kappa}
\newcommand{\A}{\mathcal{A}}
\renewcommand{\L}{\mathcal{L}}

\newcommand{\tb}{\tilde \beta}
\newcommand{\tk}{\tilde \kappa}
\newcommand{\cinc}[1]{\frac 1{2\pi\mathrm{i}}\int_{#1}}
\newcommand{\bgl}{\laxgl(n)}

\newcommand{\laxgl}{\overline{\mathfrak{gl}}}

\newcommand{\cint}{\frac 1{2\pi\mathrm{i}}\int_{C}}
\newcommand{\w}{\omega}

\newcommand{\rank}{\operatorname{rk}}
\newcommand{\rk}{\operatorname{rk}}
\newcommand{\ord}{\operatorname{ord}}

\newcommand{\tr}{\mathrm{tr}}
\newcommand{\sln}{\mathfrak{sl}}
\newcommand{\sn}{\mathfrak{s}}
\newcommand{\so}{\mathfrak{so}}
\newcommand{\spn}{\mathfrak{sp}}
\newcommand{\gl}{\mathfrak{gl}}
\newcommand{\slnb}{{\overline{\mathfrak{sl}}}}
\newcommand{\snb}{{\overline{\mathfrak{s}}}}

\newcommand{\glb}{{\overline{\mathfrak{gl}}}}
\newcommand{\fpz}{\frac {d }{dz}}
\newcommand{\lpz}{\C[z,z^{-1}]}
\newcommand{\de}{\delta}
\newcommand{\tX}{\widetilde{X}}
\newcommand{\tY}{\widetilde{Y}}

\begin{document}
\date{2010}

\title[Central Extensions of Lax operator algebra]
{Almost-graded \\ central extensions of \\ Lax operator algebras}

%%% Note the way of obtaining lower case symbols in the title %%%

\author[M. Schlichenmaier]{Martin Schlichenmaier}
\address{
University of Luxembourg, \\Mathematics Research Unit,  FSTC\\
6, rue Coudenhove-Kalergi,
L-1359 Luxembourg-Kirchberg\\
Email: martin.schlichenmaier@uni.lu }

\maketitle

\begin{abstract}

Lax operator algebras constitute a new class of infinite dimensional
Lie algebras of geometric origin. More precisely, they are
algebras of matrices whose entries are   meromorphic  functions on
a compact Riemann surface. They generalize classical
current algebras and current algebras of Krichever-Novikov type.
Lax operators for $\gl(n)$, with the spectral parameter on a 
Riemann surface,  were introduced  by Krichever. In joint
works of Krichever and Sheinman their algebraic structure was
revealed and extended to more general groups.
These algebras are almost-graded.
In this article 
their definition is recalled   and
classification and uniqueness results
for almost-graded central extensions for this new class
of algebras are presented. The explicit forms of the defining
cocycles are given.
If the finite-dimensional Lie algebra
on which  the Lax operator algebra  
is based is simple
then, up to equivalence and rescaling of the
central element, 
 there is a unique non-trivial almost-graded
central extension.
These results are  joint work
with Oleg Sheinman.

\vskip 0.5cm

{\it This is an extended  write-up of a talk presented at the 
$5^{th}$ Baltic-Nordic AGMP Workshop: Bedlewo, 12-16 October, 2009} 
\end{abstract}
%%%%%%%%%%%%%%%%%%%%%%%%%%%%%%%%%%%%%%%%%%%%%
\section{Introduction} 
%%%%%%%%%%%%%%%%%%%%%%%%%%%%%%%%%%
Classical current algebras (also called loop algebras) and their 
central extensions, the affine Lie algebras, are of fundamental 
importance in quite a number of fields in mathematics and its
applications. These algebras are examples of infinite dimensional 
Lie algebras which are still tractable. They constitute the subclass of
Kac-Moody algebras of untwisted affine type \cite{Kac}.

In the usual approach they are presented in a purely algebraic manner.
But there is a very useful geometric description behind.
In fact, the
classical current algebras correspond to Lie algebra valued 
meromorphic functions on the Riemann sphere (i.e.~on the
unique compact Riemann surface of genus zero) which are allowed to
have poles only at two fixed points.
If this rewriting is done, a very useful generalization 
(e.g.~needed in string theory, but not only there) is to consider the 
objects over a compact Riemann surface of arbitrary genus with 
possibly more than two
points where poles are allowed.
Such objects (vector fields, functions, etc.) and central extensions
for higher genus with two possible poles were 
introduced by Krichever and Novikov
\cite{KN} and generalized by me to the 
multi-point situation \cite{Schlce}.
These objects are of importance in a global operator approach to
Conformal Field Theory 
 \cite{Schlglwz}, \cite{SchlShWZ1}, \cite{SchlShWZ2}

More generally, the current algebra resp. 
their central extensions, the affine algebras, correspond to
infinite dimensional symmetries of  systems.
Moreover,  they define objects over the moduli space of Riemann surfaces
with marked points (Wess-Zumino-Novikov-Witten (WZNW) models, 
Knizhnik-Zamolodchikov (KZ) connections, etc., see e.g. \cite{SchlShWZ1}).

Their use, e.g. in quantization, 
regularizations, Fock space representation,
 force us to consider central extensions of these
algebras. 
Note that the well-known Heisenberg algebra is a 
central extension of a commutative Lie algebra.
For the classical current algebras the Kac-Moody algebras 
of affine type are obtained via central extensions.

More recently, a related class of current algebra appeared, the 
Lax operator algebras. Again 
these are  algebras of current type associated to
finite-dimensional Lie algebras and to 
compact Riemann surfaces of arbitrary genus, resp. to  smooth projective curves
over the complex numbers.
They find applications in the theory of 
 integrable systems and are connected to  the moduli space of
bundles over compact Riemann surfaces.

In 2002 Krichever \cite{Klax}, \cite{Kiso} 
studied  $\gl(n)$ Lax operators on 
higher genus Riemann surfaces. In 2007 Krichever and Sheinman
\cite{KSlax} uncovered  their algebraic structure not only for 
$\gl(n)$ but also for more general 
classes of finite-dimensional algebras, e.g. for 
$\sln(n)$, $\so(n)$ and $\spn(2n)$.

Krichever-Novikov current type algebras consist of 
Lie algebra valued meromorphic functions
on a fixed compact Riemann surface (of genus $g$) with possible poles at 
a finite set of points. The orders of the poles are 
not restricted. For the Lax operator algebra associated to 
$\gl(n)$ the elements are allowed to have additional poles
of maximal order one at a finite set of $n\cdot g$ additional 
points $\ga_s$, the weak singularities.
To each point $\ga_s$ a vector $\alpha_s\in\C^n$ is assigned 
which enters the formulation of the required form of the expansion
of the element at the point  $\ga_s$, see \refS{lax}.
The appearance of this additional poles 
is due to the fact, that the Lax operators 
operate on functions
representing sections of a non-trivial  rank $n$ vector bundle.
The additional data needed to describe the algebra are the
Tyurin parameters \cite{Tyvb} appearing in the context
of the moduli space of vector bundles over Riemann surfaces.
In this context the Krichever-Novikov current algebra
for $\gl(n)$  corresponds to
the  trivial rank $n$ vector bundle.

The classical counterparts of these algebras are graded algebras. Such
a grading is important e.g. in the context of representations,
for example to define highest weight representations.
In higher genus the Krichever-Novikov type algebras are not graded, but
only almost-graded (see \refD{almgrad}).
Fortunately the almost-gradedness is enough for the constructions
in representation theory.
It turns out that Lax operator algebras are also almost-graded.

Starting from an almost-graded Lie algebra those central extensions
are of particular interest for which the almost-grading can be
extended to the central extension.
Classification results for almost-graded extensions for the Krichever-Novikov
current algebras are given in \cite{Saff}.
For the Krichever-Novikov current algebras associated to 
finite-dimensional simple Lie algebras 
 there is up to equivalence of the extension and rescaling of the
central element only 
one nontrivial almost-graded central extension.

For the Lax operator 
algebras in joint work with Oleg Sheinman we classified 
almost-graded central extensions. 
Again it turns out that in the case that the 
associated  simple finite-dimensional
Lie algebra is
$\sln(n)$, $\so(n)$ or $\spn(2n)$
the almost-graded central extension is essentially unique.
We give its explicit form. By requiring 
a certain invariance under the action of vector fields
($\L$-invariance, see \refD{linv}) we even fix its representing two-cocycle
in its cohomology class.
The results appeared in \cite{SSlax}. It is the goal of this
contribution to report on the results.

%%%%%%%%%%%%%%%%%%%%%%%%%%%%%%%%%%%%%%%%%%%%%%%%%%%%%%%%%%%%%%%%%%
%%%%%%%%%%%%%%%%%%
\section{Krichever-Novikov type current algebras}
%%%%%%%%%%%%%%%%%%%%%%%%%%%%%%%
Let us first consider the classical situation. We fix 
 $\g$ a finite-dimensional complex Lie algebra.
The classical current algebra $\gb$
(sometimes also called loop algebra) is obtained
by tensoring $\g$ by the (associative and commutative)
algebra $\lpz$ of Laurent polynomials, i.e. 
$\gb=\g\otimes \lpz$ with the Lie product
\begin{equation*}%\label{E:classalg}
[x\otimes z^n,y\otimes z^m]:=[x,y]\otimes z^{n+m},
\quad x,y\in \g, \ n,m\in\Z.
\end{equation*}
By setting $\deg(x\otimes z^n):=n$ the Lie algebra $\gb$ becomes a 
graded algebra.
Clearly, $\gb$ is an infinite dimensional Lie algebra.
These algebras  
appear e.g. as symmetry algebras of systems with infinitely many 
independent symmetries
\cite{Kac}.
In applications quite often one is forced to consider central 
extensions of them. And we will do so also further down.

To understand in which sense the algebras of Krichever-Novikov type 
are higher genus version of both the classical current algebras and 
their central extensions we first have to geometrize the
classical situation. 
The associative algebra of Laurent polynomials $\lpz$ 
can be described  as the algebra consisting of those 
meromorphic functions on the Riemann sphere (resp. the complex
projective line $\P^1(\C)$) which are holomorphic outside
$z=0$ and $z=\infty$ ($z$ the quasi-global coordinate).
The current algebra $\gb$ can be interpreted as Lie algebra of
$\g$--valued
meromorphic functions on the Riemann sphere with possible poles
only at $z=0$ and $z=\infty$.

The Riemann sphere is the unique compact Riemann surface of
genus zero. From this point of view the next step is to take
$\Sigma$ any compact Riemann surface of arbitrary genus $g$ and
 an arbitrary finite set $A$ of points where poles 
of the meromorphic objects will be allowed.
In this way we obtain the higher genus (multi-point) current algebra
as the algebra of $\g$--valued functions on $X$ with only possible 
poles at $A$.
But we will need gradings, and later on also 
central extensions.

Some ``grading'' is essential for infinite-dimensional Lie algebras to
construct highest-weight representations, vacuum representations, etc.
In fact, for higher genus and even for genus zero  if we allow more than
two points for possible poles the  algebras under consideration will not be
graded in the usual sense. Fortunately, a weaker concept, 
an {\it almost-grading} will be enough. 
\begin{definition}\label{D:almgrad}
Let $V$ be an arbitrary Lie algebra. It is 
called {\it almost-graded} if  
\begin{enumerate}
\item
$V=\oplus_{m\in\Z}V_m$ as vector space, 
\item
$\dim V_m<\infty$, $\forall m\in\Z$,  
\item
there
exists $L_1,L_2\in\Z$ such that
\begin{equation*}
[V_k,V_m]\subseteq\bigoplus_{h=k+m-L_1}^{k+m+L_2}V_h,\qquad  \forall k,m\in \Z.
\end{equation*}
\end{enumerate}
An analogous definition works for associative algebras and for modules
over  almost-graded algebras.
\end{definition}

To introduce such a grading 
we split the set of points $A$ where poles are allowed 
into two non-empty disjoint subsets
$I$ and $O$, $A=I\cup O$.
Let $K$ be the number of points in $I$.
The points in $I$ are called {\it in-points}, the points in $O$
{\it out-points}.%
\footnote{In the interpretation of string theory, where the Riemann
  surface
$\Sigma$ corresponds to the world-sheet of the string, $I$ corresponds 
to the entry points of incoming free strings and $O$ to the 
emission points of outgoing free strings.}

%%%%%%%%%%%%%%%%%%%%%%%%%%%%%%%%%%%%%%%%%%%% 
Let $\A$ be  the associative algebra of those functions 
which are  
meromorphic on $\Sigma$ and holomorphic outside of $A$.
In some earlier work \cite{Schlce} I constructed
a certain adapted basis of $\A$
\begin{equation*}
\{A_{n,p}\mid  n\in\Z,\ p=1,\ldots,K\}.
\end{equation*}
The functions $A_{n,p}$ are uniquely fixed by giving certain vanishing orders 
at the points in $I$ and $O$ and some normalization conditions.
We will not need the exact conditions in the following. Here I only
give the vanishing order at $I$
\begin{equation*}
\ord_{P_i}(A_{n,p})=n+1-\de_{i}^p,\quad \forall P_i\in I.
\end{equation*}
Of course, a positive vanishing order means a zero, a negative a pole.

We set
$\A_n:=\langle A_{n,p}\mid p=1,\ldots, K\rangle$ for the 
$K$-dimensional subspace of $\A$. 
We call the elements of $\A_n$ homogeneous elements of degree $n$.
Clearly, 
$\A=\oplus_{n\in\Z} \A_n$ and in
 \cite{Schlce} it is shown that there exists a constant
$L$  such that
\begin{equation*}
\A_n\cdot A_m\subseteq \bigoplus_{h=n+m}^{n+m+L} \A_h,
\qquad  \forall n,m\in\Z.
\end{equation*}
Obviously, \refD{almgrad} is fulfilled and we obtain an
almost-graded structure for the algebra $\A$.

For genus zero and $I=\{0\}$, $O=\{\infty\}$
the prescription yields  $A_n=z^n$
and $\A=\lpz$. In this case the algebra is  
graded.
Note that in the two-point case we have $K=1$ and we will drop
the second subscript $p$ in $A_{n,p}$.

\begin{remark}
The notion of almost-gradedness 
was introduced by Krichever and Novikov
\cite{KN} (they called it quasi-gradedness) and they 
constructed such an almost-grading in the higher genus and 
two point case.
To find an almost-grading in the multi-point case is more difficult
and has been done in \cite{Schlth}, \cite{Schlce}.
The constant $L$  depends in a known manner on the 
genus $g$ and the number of points in $I$ and $O$.
\end{remark}
\begin{definition}
Given a finite-dimensional Lie algebra $\g$, then 
the {\it higher genus multi-point current algebra}  $\gb$ is the tensor
product
$\gb=\g\otimes \A$ with  Lie product
\begin{equation*}
[x\otimes f,y\otimes g]:=[x,y]\otimes (f\cdot g)
\end{equation*}
and almost-grading 
\begin{equation*}
\gb=\bigoplus_{n\in\Z}\gb_n, 
\qquad \gb_n:=\g\otimes \A_n.
\end{equation*}
\end{definition}
\noindent
For $g=0$ and $I=\{0\}$, $O=\{\infty\}$ this gives exactly the
classical current algebras.
See \cite{KN, Sha, Saff}.
\medskip

Additionally, we will need 
$\L$, the {\it Lie algebra of meromorphic vector fields}  on $\Sigma$ 
holomorphic outside of $A$.
This algebra is again an almost-graded algebra.
The grading is given with the help of 
a certain adapted basis 
\cite{Schlce}
\begin{equation*}
\{e_{n,p}\mid  n\in\Z,\ p=1,\ldots,K\}.
\end{equation*}
The conditions for $e_{n,p}$ are similar to the conditions
for the $A_{n,p}$. Here we only note
\begin{equation*}%\label{E:edeg}
\ord_{P_i}(e_{n,p})=n+2-\de_{i}^p,\quad \forall P_i\in I.
\end{equation*}
For genus zero and $I=\{0\}$, $O=\{\infty\}$
we get $e_{n,p}=z^{n+1}\fpz$,
and obtain in such a way as $\L$ the Witt algebra 
(sometimes also called  the Virasoro
algebra without central extension).
%%%%%%%%%%%%%%%%%%%%%%%%%%%%%%%%%%%%%%%%
%%%%%%%%%%%%%%%%%%%%%%%%%%%%%%%%%%%%%%%%

\section{Lax operator algebras}\label{S:lax} 

For higher genus there is another kind of current type algebras
given by the Lax operator algebras of higher genus.
They are related to integrable systems and to the moduli space
of semi-stable framed vector bundles over $\Sigma$.
Again 
let $\Sigma$ be a compact Riemann surface of genus $g$ 
and $A$ a  finite subset of points divided into two
nonempty disjoint subset $A=I\cup O$.
For simplicity we consider here only the {\it two-point} 
situation $I=\{P_+\}$ and $O=\{P_-\}$,
but the results are true in the more general setting
\cite{SSlaxm}.

For $n\in\N$ we fix
$n\cdot g$ additional  points on $\Sigma$ 
\begin{equation*}
W:=\{\ga_s\in\Sigma\setminus\{P_+,P_-\}\mid s=1,\ldots, n g\}.
\end{equation*}
To every point $\ga_s$ we assign a vector $\a_s\in\C^n$.
The system 
\begin{equation*}
T:=\{(\ga_s,\a_s)\in\Sigma\setminus\{P_+,P_-\}\times \C^n\mid s=1,\ldots, n g\}
\end{equation*}
is called 
a {\it Tyurin data}. 
This data is related to
the moduli of  semi-stable framed vector bundles over $\Sigma$
\cite{Tyvb}, see \refS{tyu}.

We fix local coordinates
$z_\pm$ at $P_\pm$ and $z_s$ at $\ga_s$,  $s=1,\ldots, n g$. 
Let $\g$ be  one of the matrix algebras
$\gl(n)$, $\sln(n)$, $\so(n)$, $\spn(2n)$, or $\sn(n)$, where
the latter  denotes the algebra of scalar matrices.

We consider meromorphic functions
\begin{equation}\label{E:map}
L:\ \Sigma\ \to\  \g,
\end{equation}
which are
holomorphic outside  $W\cup \{P_+, P_-\}$, have at most poles
of order one (resp. of order two for $\spn(2n)$) 
at the points in $W$, and fulfil certain 
conditions (described below) at $W$ 
depending on the Tyurin data $T$ and the Lie algebra 
$\g$.
The singularities at $W$ (resp. in abuse of notation 
the  points in $W$ themselves) are called {\it weak singularities}.

In this section we only give the conditions for the case $\g=\gl(n)$, 
$\sln(n)$, and $\sn(n)$.
The conditions for $\so(n)$ and $\spn(2n)$ are given in \refS{app}.
Let $T$ be fixed.
For $s=1,\ldots, ng$ we require that there exist $\b_s\in\C^n$ 
and $\ka_s\in \C$ such that the
function $L$ has an expansion at $\ga_s\in W$ of the form 
\begin{equation}\label{E:glexp}
L(z_s)_{|}=\frac {L_{s,-1}}{z_s}+
L_{s,0}+\sum_{k>0}L_{s,k}{z_s^k}
\end{equation}
with
\begin{equation}\label{E:gldef}
L_{s,-1}=\a_s \b_s^{t},\quad
\tr(L_{s,-1})=\b_s^t \a_s=0,
\quad
L_{s,0}\,\a_s=\ka_s\a_s.
\end{equation}
In particular, 
if $\a_s$, and $\b_s\ne 0$ the matrix $L_{s,-1}$ is a rank 1 matrix, and  
$\a_s$ is  
an eigenvector of $L_{s,0}$.
In \cite{KSlax} it is shown that
the requirements \refE{gldef} are independent of the chosen
coordinates $z_s$ and that 
the set of all such functions constitute an associative algebra under
the point-wise matrix multiplication. 
We denote it by $\gb$.%
\footnote{With the intent not too overload the notation we
use $\gb$ for any of the current algebra versions associated to
$\g$. It will be clear from the context whether it is a classical,
a Krichever-Novikov, or a Lax operator algebra.}
As Krichever and Sheinman \cite{KSlax} showed $\gb$ will always be a
Lie algebra under the matrix commutator.
The main point is to verify that the pole orders at the weak singularities 
do not increase and that the expansion is of the required type.
(Note that the coefficients in the expansion \refE{glexp}
are matrices and the conditions \refE{gldef} 
have  to be used in the verification.)
%%%%%%%%%%%%%%%%%%%%%%%%%%%%%%%%%
\begin{remark}
If $\a_s=0$ for all $s$ there are no additional singularities and
we obtain the usual Krichever-Novikov current algebras back. 
\end{remark}
\begin{remark}
For the subalgebra $\sln(n)$ of $\gl(n)$ by \refE{map} all matrices $L_{s,k}$
in \refE{glexp} have to be traceless. The conditions 
\refE{gldef} stay the same.
\end{remark}
\begin{remark}
For the subalgebra $\sn(n)$ 
all matrices have to be scalar matrices. As $L_{s,-1}$ is 
traceless hence,  either $\a_s=0$ or 
$\b_s=0$. In both cases there is no pole at $\ga_s$.
Furthermore the eigenvalue condition for $L_{s,0}$ is 
also true.
Hence  $\snb(n)$  coincides with the Krichever-Novikov function
algebra, i.e. 
\begin{equation*}%\label{E:sint}
\snb(n)\cong \sn(n)\otimes \A\cong \A,
\end{equation*}
as associative algebras.
\end{remark} 
In fact we have a splitting
$\gl(n)=\sn(n)\oplus \sln(n)$ given by
\begin{equation*}
X\mapsto \left(\ \frac {\tr(X)}{n}I_n\ ,\ X-\frac {\tr(X)}{n}I_n\ \right),
\end{equation*}
where $I_n$ is the $n\times n$-unit matrix.
This splitting 
induces a corresponding splitting for the  Lax operator 
algebra  $\glb(n)$:
\begin{equation*}
 \glb(n)=\snb(n)\oplus \slnb(n).
\end{equation*}

\medskip
\subsection{The almost-grading}
In contrast to the Krichever-Novikov situation the almost-grading
of $\A$ cannot be directly used to 
introduce an almost-grading  for the Lax operator algebras (with
$\a_s\ne 0$).
But similar ideas for introducing the almost-grading in $\A$ and $\L$
work here too.
For every $m\in\Z$  a subspace $\gb_m$ inside  $\gb$  is
defined  as the subspace where non-zero elements 
have matrix expansions of order $m$ at $P_+$, and
a corresponding order at $P_-$ forcing the elements  to be
essentially unique.
In \cite{KSlax} it is shown that
$\gb$ is almost-graded. More precisely, we have
\begin{equation*}
\gb=\bigoplus_{m\in\Z}\gb_m,\quad  \dim\gb_m=\dim \g,\quad
[\gb_m,\gb_k]\subseteq\bigoplus_{h=m+k}^{m+k+M}\gb_h,
\end{equation*}
with a constant $M$ independent of $m$ and $k$. In fact, for generic
situations $M=g$, the genus of $\Sigma$, will work.
To give an idea, for generic $m$ 
the element $L\in\gb$, $L\ne 0$ lies in $\gb_m$ if there exists 
$X_+,X_-\in\g\setminus\{0\}$ with 
\begin{equation}\label{E:degree}
L(z_+)_|=X_+z_+^m+O(z_+^{m+1}),\qquad
L(z_-)_|=X_-z_-^{-m-g}+O(z_-^{-m-g+1}).
\end{equation}
As $\dim\gb_m=\dim\g$ we obtain that given $X\in\g$  
there exists a unique $X_m\in\gb_m\subset\gb$ such that
(see also \cite[Prop.~2.4]{SSlax})
\begin{equation}\label{E:degree1}
X_m(z_+)_|=X\cdot z_+^m+O(z_+^{m+1}).
\end{equation}

%%%%%%%%%%%%%%%%%%%%%%%%%%%%%%%%%%%%%%%%%%%%%%%%%%%%%%
\section{The geometric meaning of  Tyurin parameters}
\label{S:tyu}
%%%%%%%%%%%%%%%%%%%%%%%%%%%%%%%%%%%%%%%%%%%%%%%%%%%
Despite the fact that in this article 
I will not use the geometric relevance of 
Tyurin parameters in relation to the moduli space of bundles,
it might be interesting to recall this background information.
The reader not interested in this connection might directly jump to
the next section.

Let $\Sigma$ be a compact Riemann surface (or in the language of
algebraic geometry a projective smooth curve over $\C$)
of genus $g$.
Fix a number $n\in\N$. Given a rank $n$ holomorphic (resp. algebraic)
vector bundle $E$ its determinant $\det E$ 
is defined as $\det E=\wedge^nE$.
The degree $\deg E$ of the bundle $E$ is defined as 
$\deg(\det(E))$. Recall that for a line bundle $M$ 
over a Riemann surface the degree of
$M$ can be determined by taking a global meromorphic section of $M$
and counting the number of zeros minus the number of poles of this 
section.

For vector bundles over compact Riemann surfaces  we have the 
Hirzebruch-Riemann-Roch formula
\begin{equation}\label{E:HRR}
\dim \mathrm{H}^0(\Sigma,E)-
\dim \mathrm{H}^1(\Sigma,E)=
\deg E-\rank(E)(g-1).
\end{equation}

If one wants to construct a moduli space for vector
bundles
one has to restrict the set of  vector bundles to the subset
of stable, or more general semi-stable bundles.
\begin{definition}
A bundle $E$ over a projective smooth curve is 
called {\it stable}, if for all subbundles $F\ne E$ one has
\begin{equation}\label{E:stab}
\frac {\deg F}{\rk F}\quad <\quad \frac {\deg E}{\rk E}.
\end{equation}
The bundle $E$ is called {\it semi-stable} if 
in \refE{stab} the strict inequality $<$ is replaced by $\le$.
\end{definition}

In the following we consider bundles $E$ which are of rank $n$ and 
degree $n\cdot g$.
If we evaluate \refE{HRR} for such bundles we 
obtain 
the value $n$ on the left hand side.
For a generic semi-stable bundle $E$ one has 
$\dim \mathrm{H}^1(\Sigma,E)=0$, hence we get
\begin{equation*}
\dim \mathrm{H}^0(\Sigma,E)=n.
\end{equation*}
If we choose a basis $S:=\{s_1,s_2,\ldots, s_n\}$ of the 
space of global holomorphic sections
of $E$, their exterior power is a global holomorphic section 
\begin{equation*}
s_1\wedge s_2\cdots\wedge s_n\ \in\  \mathrm{H}^0(\Sigma,\det E).
\end{equation*}
The zeros  of this section are exactly the points $P\in\Sigma$
for which  the set  $S$ fails to be a basis of the
fibre 
of $E$ at the point $P$. This fibre is denoted by
$E_P$.
As the degree  $\deg E=ng$ there exist (counted with 
multiplicities) exactly $ng$ such points. 
For a generic choice of the set of sections, all zeros will be
simple. Hence, we will obtain  $ng$ such points.
They correspond exactly to the weak singularities $W$ 
appearing in the definition of the Lax operator algebras.
Accordingly we denote the zero points by $\ga_s$, 
$s=1,\ldots,ng$. 

Furthermore, as the zero is of order one at such a $\ga_s$
the sections evaluated at $\ga_s$ 
span a $(n-1)$-dimensional subspace of $E_{\ga_s}$.
We have relations
\begin{equation*}
\sum_{i=1}^{n} \alpha_{s,i} s_i(\ga_s)=0,\qquad s=1,\ldots, ng.
\end{equation*}
In this way to every $\ga_s$ a vector $\alpha_s\in\C^n$,
$\alpha_s\ne 0$  can be assigned.
This vector is unique up to multiplication by a non-zero scalar, hence
unique as element $[\alpha_s]\in\P^{n-1}(\C)$.
Again these vectors are exactly the vectors  used in the definition
of the Lax operator algebras.
If one checks the conditions \refE{gldef} one sees immediately that
they are independent 
of a rescaling. Hence only the 
projective class $[\alpha_s]$ matters.

Obviously everything depends on the set $S$ of basis elements.
The choice of such a basis is called a {\it framing} of the bundle $E$.
In the way  described above the space of Tyurin parameters 
  parameterizes an open dense subset of semi-stable framed 
vector bundles of rank $n$ and degree $ng$.
Note also that given such a bundle $E$ with  
fixed set $S$ of basis elements of $\mathrm{H}^0(\Sigma,E)$  
it can be trivialized over $\Sigma\setminus W$.

Associated to this moduli spaces integrable hierarchies 
of Lax equations can be constructed.  
See \cite{Klax} and \cite{Shih} for results in this
directions.

%%%%%%%%%%%%%%%%%%%%%%%%%%%%%%%%%%%%%%%%%%%%%%
\section{Central extensions}
%%%%%%%%%%%%%%%%%%%%%%%%%%%%%%%%%%
By the applications in quantum theory (but not only there)
we are forced to consider central extensions of
the introduced algebras. An example how they appear is given by
regularization of a not well-defined action. The regularization
makes the action well-defined but now it is only a projective
action. To obtain a linear action we have to pass to a 
suitable central extension of the Lie algebra.

In the following definition $\gb$ could be any Lie algebra of current
type or even more general any almost-graded Lie algebra.
Central extensions are given by  Lie algebra 2-cocycles
with values in the trivial  module $\C$.
Recall that such a
2-cocycle for $\gb$ is a bilinear form
$\ga:\gb\times\gb\to\C$ which is (1) antisymmetric
and (2) fulfils the cocycle condition
\begin{equation*}%\label{E:cohcoc}
\ga([L,L'],L'')+
\ga([L',L''],L)+
\ga([L'',L],L')=0, \qquad \forall L,L',L''\in\gb.
\end{equation*}
A 2-cocycle $\ga$ is a coboundary if there exists a linear form $\phi$ 
on $\laxg$ with 
\begin{equation*}
\ga(L,L')=\phi([L,L']),\qquad  \forall L,L'\in\laxg.
\end{equation*}
Given two  cocycles $\ga$ and $\ga'$ then they are cohomologous if
their difference is a coboundary.

Given a 2-cocycle $\ga$ for $\laxg$, the associated central extension
$\gh_{\ga}$ is given as vector space direct sum $\gh_\ga=\gb\oplus\C$
with Lie product given by
\begin{equation*}%\label{E:centextf}
[\widehat{L},\widehat{L'}]=\widehat{[L,L']}+\ga(L,L')\cdot t,
\quad [\widehat{L},t]=0,\qquad
L,L'\in\laxg.
\end{equation*}
Here we used $\widehat{L}:=(L,0)$ and $t:=(0,1)$.
Vice versa, every central extension  
\begin{equation*}
\begin{CD}
0@>>>\C@>i_2>>\gh@>p_1>>\gb@>>>0,
\end{CD}
\end{equation*}
defines  a 2-cocycle
$\ga:\gb\to\C$ after choosing a (linear) section $s:\gb\to \gh$ of $p_1$
by the condition
\begin{equation*}
[s(L),s(L')]-s([L,L'])=(0,\ga(L,L')).
\end{equation*}
Different sections $s_1$ and $s_2$ give cohomologous 2-cocycles 
$\ga_1$ and $\ga_2$.
Two central extensions $\gh_{\ga}$ and $\gh_{\ga'}$ are equivalent
if the defining cocycles $\ga$ and $\ga'$ are cohomologous.

Coming from the applications 
we are only interested in those central extensions 
which allow the extension of the almost-grading of 
$\gb$ to the central extension $\gh_\ga$ by assigning to the central
element $t$ a certain degree and using the degree for $\gb$ for 
the subspace $(\gb,0)$ in $\gh$.
This is  possible if and only if our defining cocycle $\ga$ is {local}
in the following sense:
\begin{definition}
A cocycle $\ga$ for the almost-graded Lie algebra $\gb$ is {\it local}
 if and
only if there exists constants $M_1,M_2$ such that
\begin{equation*}
\ga(\gb_m,\gb_k)\ne 0\implies M_1\le m+k\le M_2,\qquad\forall m,k\in\Z.
\end{equation*}
A central extension obtained by a local cocycle and 
with the extended grading is  called an {\it almost-graded central
extension} of $\gb$. 
\end{definition}

The question is  how to construct cocycles?
For the current algebras we first fix  
$\langle .,.\rangle$ 
an invariant symmetric bilinear form on 
the finite-dimensional Lie algebra $\g$.
Recall that invariance means that 
$$\langle[a,b],c\rangle=\langle a,[b,c]\rangle
\quad \text{ for all} \quad  a,b,c\in\g.$$
For a simple Lie algebra the Cartan-Killing form is 
the only such form  up to 
rescaling. Moreover it is non-degenerate. 
For $\sln(n)$ the Cartan-Killing form is given by
$ \alpha(A,B)=\tr(AB)$.

For the classical current algebras
associated to a simple finite-dimensional Lie algebra with
Cartan-Killing form  $\langle .,.\rangle$ 
a non-trivial central extension $\gh=\gb\oplus\C t$
is defined by
\begin{equation*}%\label{E:centralclass}
[x\otimes z^n,y\otimes z^m]=[x,y]\otimes z^{n+m}-\langle x,y\rangle
\cdot n\cdot \de_{n}^{-m}\cdot t.
\end{equation*}
To avoid cumbersome notation I dropped the $\ \hat{ }\  $ here.
It is called the (classical) affine Lie algebra associated to $\g$.
Another name is untwisted affine Lie algebra of Kac-Moody type
\cite{Kac}. 
By setting $\deg t:=0$ (and using $n=\deg(x\otimes z^n)$) the 
affine algebra is a  graded algebra.
The cocycle 
$$
\ga(x\otimes z^n,y\otimes z^m)=
 -\langle x,y\rangle
\cdot n\cdot \de_{n}^{-m}
$$ is obviously local.
Indeed for $\g$ simple it is  the
only non-trivial extension up to equivalence and rescaling. 

\medskip

For the usual Krichever-Novikov algebras $\g\otimes \A$ a 2-cocycle
can be given by
\begin{equation*}
\ga_C(x\otimes f,y\otimes g):=\langle x,y\rangle \cint fdg,
\end{equation*}
where   $C$ is a closed curve on $\Sigma$ \cite{Saff}.

The cocycles depend crucially on the choice of the 
integration path $C$ and we might obtain different non-cohomologous
cocycles, hence non-equivalent central extensions.
But recall that we are mainly interested in local cocycles.
There is a special class of integration cycles, the so called
separating cycles $C_s$, 
which separate the points in $I$ from the points in $O$
with multiplicity one.
As all separating cycles are homologous
the value of the integration does not depend on the  separating cycle we take.
In particular, we could take circles
around the points in $I$,
hence calculate the integral using 
residues there.
In \cite[Theorem 4.6]{Saff} 
a complete classification of local cohomology classes, i.e.
classes which admit at least one representing element which is local,
is given for the case of reductive Lie algebras $\g$. 
In this way a classification of almost-graded central extensions
of the Krichever-Novikov current algebras is obtained.
In particular, for $\g$ simple there is 
only one non-trivial almost-graded extension
 up to equivalence and
rescaling.
Let me stress the fact that without the  condition of locality of the
cocycle the statement would be wrong.
%%%%%%%%%%%%%%%%%%%%%%%%%%%%%%%%%%%%%%%%%%%%%%%%%%%%%%%%
%%%%%%%%%%%%%%%%%%%%%%%%%%%%%%%%
\section{Central extensions for Lax operator algebras}
%%%%%%%%%%%%%%%%%%%%%%%%%%%%%%%%%%%%%%%%%%%%%%%%%%%
For the Lax operator algebras we obviously 
have the problem that differentiation  of our objects
will increase the pole order at the weak singularities. We
will not stay in the algebra. 
The deeper reason for this is that the objects are not
really functions but representing functions of sections 
of a bundle.
To correct this we  need
a connection $\nabla$ and have to take the covariant derivative.
It will be defined with the help of a connection form $\w$.
This form is  a $\g$-valued 
meromorphic 1-form, holomorphic outside $P_+$, $P_-$ and $W$, and
with prescribed behavior at the points in $W$.
For $\gamma_s\in W$ with $\a_s= 0$ the requirement is that
$\w$ is also regular there.
For the points $\gamma_s$ with   $\a_s\ne 0$ we require that
it has  the  expansion
\begin{equation*}%\label{E:connl}
\w(z_s)_|=\left(\frac {\w_{s,-1}}{z_s}+\w_{s,0}+
\sum_{k>1}\w_{s,k}z_s^k\right)dz_s.
\end{equation*}
For $\gl(n)$ we require that there
exist $\tb_s\in\C^n$ 
and $\tk_s\in \C$ such that 
\begin{equation*}%\label{E:gldefc}
\w_{s,-1}=\a_s \tb_s^{t},\quad
\w_{s,0}\,\a_s=\tk_s\a_s,
\quad
\tr(\w_{s,-1})=\tb_s^t \a_s=1.
\end{equation*}
Note that compared to \refE{gldef} only the last condition was
modified.
These conditions were introduced in 
\cite{Klax}, \cite{KSlax}, see also \cite{SSlax}. 
We will even require that $\w$ is holomorphic at the point $P_+$
(resp.
at the points in $I$). As we allow poles at the point $P_-$
(resp. at the points in $O$) such a $\w$ will always exist.

The induced  connection, resp. covariant derivative for the algebra will be
given by 
\begin{equation*}%\label{E:conng}
\nabla^{(\w)}=d+[\w,.],\qquad
\nabla_e^{(\w)}=dz(e)\frac {d}{dz}+[\w(e),.\,].
\end{equation*}
Here $e$ is a vector field from $\L$.

\begin{remark}
For the subalgebras $\slnb(n)$ 
and $\snb(n)$
we can take the same $\w$ as for $\glb(n)$. The covariant derivative
will respect the subalgebras. In fact for the scalar algebra 
$\nabla^{(\w)}=d$.
For the other
algebras see \refS{app}.
\end{remark}

\begin{proposition}\cite{SSlax}\label{P:deralm}
The covariant derivative $\nabla_e^{(\w)}$ acts 
as a derivation on $\gb$ and  makes 
$\gb$  to  an almost-graded Lie module over  $\L$.
\end{proposition}

In \cite{SSlax} the following 2-cocycles for $\gb$ were given
\begin{equation}\label{E:g1}
\ga_{1,\w,C}(L,L'):= \cinc{C} \tr(L\cdot \nabla^{(\w)}L'),
\qquad L,L'\in\laxg, 
\end{equation}
and
\begin{equation}\label{E:g2}
\ga_{2,\w,C}(L,L'):= \cinc{C} \tr(L)\cdot \tr(\nabla^{(\w)}L'),
\qquad L,L'\in\laxg.
\end{equation}
Here $C$ is an arbitrary closed path in $\Sigma$.
Indeed  these are cocycles. The cocycle
$\ga_{2,\w,C}$ does not depend on $\w$ and will vanish for
$\g\ne\gl(n),\sn(n)$.
Two cocycles  $\ga_{1,\w,C}$ and   $\ga_{1,\w',C}$ with different connection
forms $\w$ and $\w'$ will be cohomologous.
\begin{definition}\label{D:linv}
Let the action of $\L$ on $\gb$ be given via $\nabla^{(\w)}$.
A cocycle $\ga$ is called $\L$-invariant if and only if
\begin{equation}\label{E:linv}
\ga(\nabla_e^{(\w)}L,L')+
\ga(L,\nabla_e^{(\w)}L')=0,\qquad\forall e\in\L.
\end{equation}
\end{definition}
\noindent
The  2-cocycles \refE{g1},\refE{g2} are
$\L$-invariant.

\medskip
The cocycles depend crucially on the choice of the 
integration path $C$.
As in the Krichever-Novikov case we are interested only in 
local cocycles.
Hence, we take again the 
separating cycles $C_s$ as $C$ .
The integral does not depend on the separating cycle we take
(see \cite[Prop. 3.6]{SSlax}). 
Note as  here we have additional poles ``between $I$ and $O$''
there is indeed something to prove.
The integral over $C_s$ can be determined by calculating the 
residue at $P_+$ (resp. at $P_-$). Indeed, for $C=C_s$ the cocycle
\refE{g1} in  a modified form can already be found in 
\cite{KSlax}.

If we take a separating cycle  as integration path we drop the
reference to $C_s$ for the cocycle, i.e. we use $\ga_{1,\w}$ and $\ga_2$.
\begin{proposition}\cite{SSlax}
The cocycles 
$\ga_{1,\w}$, and $\ga_2$ are $\L$-invariant local cocycles.
\end{proposition}
But what about the opposite, i.e. is every local and $\L$-invariant
cocycle a linear combination of these two cocycles? One of the main
results of \cite{SSlax} is that this is true.
To formulate the results we need to introduce the following convention.
A cohomology class will be called local, resp. $\L$-invariant if it 
admits a representing cocycle which is  local, resp. $\L$-invariant.
This does not imply that all representing cocycles will be of this 
type.
\begin{theorem}\label{T:main}
\cite[Thm. 3.8]{SSlax}
$ $

(a) If $\g$ is simple (i.e. $\g=\sln(n), \so(n), \spn(2n)$) then 
the space of local cohomology  classes for $\gb$ is one-dimensional.
If we fix any connection form $\w$ then the space 
will be generated by the class of $\ga_{1,\omega}$.
Every $\L$-invariant (with respect to the connection 
$\w$) local cocycle is a scalar multiple of 
$\gamma_{1,\omega}$.

(b) For $\laxg=\bgl$ the space of local 
cohomology classes which are $\L$-invariant
having been restricted to the scalar subalgebra is two-dimensional.
If we fix any connection form $\w$ then the space
will  be generated by the classes of the 
cocycles $\gamma_{1,\omega}$ and $\gamma_2$.
Every $\L$-invariant local cocycle is a linear combination of  
$\gamma_{1,\omega}$ and $\gamma_2$.
\end{theorem}

\begin{corollary}
Let $\g$ be a simple classical Lie algebra and $\gb$ the associated
Lax operator  algebra. 
Let $\w$ be a fixed connection form.
Then in each local cohomology class
$[\ga]$ there exists a unique 
representative $\ga'$ which is local and $\L$-invariant (with respect to
$\w$). 
Moreover, $\ga'=\a\ga_{1,\omega}$, with $\a\in\C$.
\end{corollary}
\begin{remark}
By this corollary for $\gb$ with $\g$ simple every local cohomology class
will be $\L$-invariant. This is not true for $\glb(n)$, due to its
abelian part $\snb(n)$.
For  $\snb(n)\cong\A$ the coboundaries are zero and the
cocycle condition reduces to the antisymmetry of the
bilinear form.  Hence it is possible to write down a lot of
different  local cocycles which are not equivalent. 
Only by     the $\L$-invariance the cocycle will be uniquely fixed
up to multiplication by a scalar.
As explained in \cite{Scocyc} $\L$-invariance is a natural condition 
as it is automatically true for those representations of
$\A$ which are in fact representations of the larger algebra
of differential operators of degree $\le 1$.
\end{remark}
%%%%%%%%%%%%%%%%%%%%%%
\begin{remark}
The notion of $\L$-invariance should be compared to the notion
of $S^1$-invariant cocycles for the classical loop algebras
\cite{PS}.
\end{remark}
%%%%%%%%%%%%%%%%%%%%%%%%%%%%%%%%%%%%%%%%%%%%%%%%%%%%%
 
\section{Some ideas of  the proof}
%%%%%%%%%%%%%%%%%%%%%%%%%%%%%%%%%%%%%%%%%%%%%%%%%%%%%%%%%
The proofs are technically quite involved
and can be found in \cite{SSlax}. 
Nevertheless  I want to
give some of the principal ideas.

\medskip

{\bf (a)} 
We start with a local, $\L$ invariant cocycle and use the almost-graded
$\L$-module structure to show that everything can be reduced to level 
zero. 
This should be understood in the following sense.
Let
$L_m$ and $L_k'$ be homogenous elements of degree $m$ and $k$
respectively, then  the level of the pair $(L_m,L_k')$ is the
sum of their degrees $m+k$. 
We show that if the level $l=m+k\ne 0$ the cocycle
values for pairs of elements of level $l$ can be linearly expressed by
values of the cocycle evaluated for pairs at higher level
with universal 
coefficients only depending on the algebra $\gb$.
By locality the cocycle values are zero for high enough level.
Hence, they will vanish for all levels $l>0$.
Moreover, by recursion their values at level $l<0$ are fixed by knowing the
values at level zero.
See below some more details on this step.

\medskip

{\bf (b)}
Next we show that for the level zero the cocycle under consideration will
be  a linear combination of the cocycles $\ga_{1,\w}$ and $\ga_2$.
Hence they will coincide everywhere.
This even gives an identification on the cocycle level, not only on
the cohomology class level.

\medskip

{\bf (c)}
The abelian part of the algebra is now covered, as we put $\L$-invariance
into the requirement.
For the simple part, we have to show that in every 
cohomology class there is an
$\L$-invariant cocycle.
To show this we consider the Chevalley generators 
and the Chevalley-Serre relations of the
finite-dimensional simple Lie algebra $\g$ and use the
almost-graded structure inside $\gb$ and the boundedness from above
of the cocycle.
We make appropriate cohomological changes and end up with the fact
that everything depends linearly on only one single 
cocycle value at level zero.
In fact everything is uniquely fixed with respect to the
value of $\ga(H_1^\alpha,H_{-1}^\alpha)$ where $\alpha$ is one 
(arbitrary) fixed simple root and $H^\alpha$ the corresponding 
generator in the Cartan subalgebra.
Hence the space is  at most one-dimensional.
But $\ga_{1,\w}$ is a local 
 and $\L$-invariant cocycle which is not a boundary
and we get the result.
As a side effect we obtain that in every local cohomology class
(for $\g$ simple) there is a unique $\L$-invariant representing cocycle.

%%%%%%%%%%%%%%%%%%%%%%%%%%%%%%%%%%%%%%%%%%%%%%%%%%

\medskip

In the following I will indicate for the interested
reader some {\bf  more facts about part (a)} above.
In particular I want to 
show the
usage of the almost-grading.

Let $\{X^r\mid r=1,\ldots,\dim\g\}$ be a basis of the
finite-dimensional
Lie algebra. As mentioned earlier we can 
find elements   
$\{X_n^r\mid r=1,\ldots,\dim\g\}$ in $\gb$
which are of order $n$ at the  point $P_+$ and start with 
leading matrix coefficient $X^r$, see \refE{degree1}.
We have the decomposition $\gb=\oplus_{n\in\Z}\gb_n$ into
subspaces of homogenous elements of  degree $n$ where 
the subspace $\gb_n$ is generated by the basis
$\{X_n^r\mid r=1,\ldots,\dim\g\}$.
The vector field $e_n\in \L$  of degree $n$ starts with 
order $n+1$ at the point $P_+$.
By \refP{deralm} the algebra  
$\gb$ is an almost-graded module over $\L$. 
Calculations in local coordinates yield 
\begin{equation}\label{E:covalm}
\nabla_{e_k}^{(\w)}X_n^r=nX_{n+k}^r+L,
\end{equation}
where $L$ is an element of $\gb$ of leading order $\ge n+k+1$
at $P_+$.
Recall that we have chosen our connection form to be holomorphic
at the point $P_+$.

For a cocycle $\ga$ evaluated for pairs of elements of level
$l$ we will use the symbol $\equiv$ to denote that the expressions are
the same on both sides of an equation up 
to universal expressions in the 
values of $\ga$  at
higher level. 
This has to be understood in the following strong sense:
\begin{equation*}
\sum \alpha^{n}_{r,s}\ga(X_n^r,X_{l-n}^s)\equiv 0,\qquad 
 \alpha^{n}_{r,s}\in\C
\end{equation*}
means a congruence modulo a linear combination of values of $\ga$ at 
pairs of basis elements of level $l'>l$. The coefficients of that
linear combination, as well as the  $\alpha^{n}_{r,s}$, depend only  
on the structure of the  Lie algebra $\gb$ and do not depend on $\ga$.

\medskip

By the $\L$-invariance \refE{linv} we have 
\begin{equation*}
\ga(\nabla_{e_p}^{(\w)}X_m^r,X_n^s)+
\ga(X_m^r,\nabla_{e_p}^{(\w)}X_n^s)\ =\ 0.
\end{equation*}
Using the almost-graded structure \refE{covalm}
we obtain the  formula 
\begin{equation}\label{E:recform}
m\ga(X_{p+m}^r,X_n^{s})+
n\ga(X_{m}^r,X_{n+p}^{s})\equiv 0,
\end{equation}
valid for all $n,m,p\in\Z$.

If  we set $p=0$ in \refE{recform} then we obtain
\begin{equation}\label{E:recform0}
(m+n)\ga(X_{m}^r,X_n^{s})
\equiv 0.
\end{equation}
Hence, for the level $(m+n)\ne 0$ everything is 
determined by the values  at higher level.
This implies in particular that if a cocycle is bounded from above it
will be automatically bounded by zero and if it vanishes at level 
zero it will vanish identically.

By evaluating \refE{recform} for suitable values for $m,p,k$ and
using the fact that all values in level greater than zero vanishes we 
obtain
\begin{equation*}%\label{E:zerodg}
\ga(X^r_m,X_0^s)= 0,\qquad
\forall  m\ge 0.
\end{equation*}
\begin{equation*}%\label{E:nm1}
\ga(X_{n}^r,X_{-n}^s)=n\cdot \ga(X_{1}^r,X_{-1}^s),
\qquad
\forall  m\in\Z.
\end{equation*}
\begin{equation*}%\label{E:p1m1}
\ga(X_{1}^r,X_{-1}^s)=\ga(X_{1}^s,X_{-1}^r).
\end{equation*}
Hence everything depends only on the 
values of $\ga(X_{1}^r,X_{-1}^s)$, with $r,s=1,\ldots,\dim\g$.

Let $X\in\g$ then we denote by 
$\tX_n$ any element
in $\gb$ with leading term $Xz_+^n$ 
at $P_+$.
We define 
\begin{equation}\label{E:cocar}
\psi_{\ga}:\g\times\g\to\C\qquad
\psi_{\ga}(X,Y):=\ga(\tX_1,\tY_{-1}).
\end{equation}
As the cocycle vanishes for level greater zero, $\psi_\ga$ does not depend
on the choice of  $\tX_1$ and $\tY_{-1}$.
Obviously, it is a bilinear form on $\g$.
Moreover in \cite{SSlax} we show
that it is symmetric and invariant.

As the cocycle $\ga$ is fixed by the 
values $\ga(X_{1}^r,X_{-1}^s)$, and they are fixed by the bilinear
map $\psi_{\ga}$ we obtain:
\begin{theorem}\cite{SSlax}
Let $\gamma$ be an $\L$-invariant cocycle for $\gb$ which is bounded 
from above by zero. Then $\ga$ is completely fixed by 
the  associated symmetric and invariant
bilinear form $\psi_\ga$ on $\g$ defined via \refE{cocar}.
\end{theorem}

In the case of a simple Lie algebra we are done, as every such   
form is a multiple of the Cartan-Killing form. And we obtain the
uniqueness (up to multiplication with a scalar) of a local 
and $\L$-invariant cocycle.
For $\gl(n)$ we  use the splitting into  $\sn(n)\oplus\sln(n)$
and have to refer to uniqueness results for the scalar algebra
$\A$ obtained in \cite{Scocyc}.

%%%%%%%%%%%%%%%%%%%%%%%%%%%%%%
\section{Appendix: $\so(n)$ and $\spn(2n)$}
%%%%%%%%%%%%%%%%%%%%%%%%%%%%%%%%%%%%%%%%%%%%%%%%%%%%%%%
\label{S:app}
\subsection{$\so(n)$}
%%%%%%%%%%%%%%%%%%%%%%%%%%%%%

In the case of {\bf $\so(n)$} we require that
all $L_{s,k}$ in \refE{glexp} are  skew-symmetric.
In particular, they are trace-less.
The set-up has to be slightly modified following \cite{KSlax}.
First only  those Tyurin parameters $\a_s$ are allowed which satisfy
$\a_s^t\a_s=0$.
Then, \refE{gldef} is modified in  the following way:
\begin{equation}\label{E:sodef}
L_{s,-1}=\a_s\b_s^t-\b_s\a_s^t,
\quad
\tr(L_{s,-1})=\b_s^t\a_s=0,
\quad
L_{s,0}\,\a_s=\ka_s\a_s.
\end{equation}
The relations \refE{sodef} do not depend on the coordinates $z_s$.
Under the point-wise
matrix commutator the set of such maps constitute a Lie algebra,
see  \cite{KSlax}.

As far as the connection form 
for $\so(n)$ is concerned 
we require that there 
exist $\tb_s\in\C^n$ 
and $\tk_s\in \C$ such that
\begin{equation*}%\label{E:sodefc}
\w_{s,-1}=\a_s\tb_s^t-\tb_s\a_s^t,
\quad
\w_{s,0}\,\a_s=\tilde\ka_s\a_s,
\quad
\tb_s^t\a_s=1.
\end{equation*}
Such  connection
forms exist.

\subsection{$\spn(2n)$}
%%%%%%%%%%%%%%%%%%%%%%%%
For {\bf $\spn(2n)$}
we consider  a symplectic form  $\hat\sigma$  
for $\C^{2n}$ given by
a non-degenerate skew-symmetric matrix $\sigma$.
Without loss of generality we might even assume 
that this matrix is given in 
the 
standard form
$\sigma=\begin{pmatrix} 0& I_n
\\ -I_n&0
\end{pmatrix}
$.
The Lie algebra $\spn(2n)$ is the Lie algebra of
matrices $X$ such that $X^t\sigma+\sigma X=0$.
This is equivalent to $X^t=-\sigma X\sigma^{-1}$, which implies that
$\tr(X)=0$.
For the standard form above, 
$X\in\spn(2n) $ if and only if
\begin{equation*}
X=\begin{pmatrix} A&B
\\
C&-A^t
\end{pmatrix}, \qquad B^t=B,\quad C^t=C.
\end{equation*}
At the weak singularities we have the expansion with matrices 
$L_{k,s}\in\spn(2n)$
\begin{equation*}%\label{E:glexpsp}
L(z_s)=\frac {L_{s,-2}}{z_s^2}+\frac {L_{s,-1}}{z_s}+
L_{s,0}+L_{s,1}{z_s}+\sum_{k>1}L_{s,k}{z_s^k}.
\end{equation*}
The conditions \refE{gldef} are  modified as
follows (see \cite{KSlax}):
there exist $\b_s\in\C^{2n}$, 
$\nu_s,\ka_s\in\C$ such that 
\begin{equation}\label{E:spdef1}
L_{s,-2}=\nu_s \a_s\a_s^t\sigma,\quad
L_{s,-1}=(\a_s\b_s^t+\b_s\a_s^t)\sigma,
\quad{\b_s}^t\sigma\a_s=0,\quad
L_{s,0}\,\a_s=\kappa_s\a_s,
\end{equation}
and 
\begin{equation}\label{E:spdef2}
\a_s^t\sigma L_{s,1}\a_s=0.
\end{equation}
Again in \cite{KSlax} it is shown that under the point-wise
matrix commutator the set of such maps constitute a Lie algebra.

For the connection form for 
$\spn(2n)$ we require  that 
there exists 
$\tb_s\in\C^{2n}$, 
$\tilde\ka_s\in\C$ such that 
\begin{equation*}%\label{E:spdefc}
\w_{s,-1}=(\a_s\tb_s^t+\tb_s\a_s^t)\sigma,
\quad
\w_{s,0}\,\a_s=\tilde\kappa_s\a_s,
\quad 
\a^t_s\sigma\w_{s,1}\a_s=0,\quad
\tb_s^t\sigma\a_s=1.
\end{equation*}

\begin{remark}
The reader might ask why for $\spn(2n)$ 
there appear poles of order two at the 
weak singularities. By direct calculations it turned out that such a
modification
has to be done to retain  that 
their will be no additional 
degree of freedom
(compared to other points)  at each weak singularity individually
(due to the  relations 
\refE{spdef1} and \refE{spdef2}), 
without disturbing the
closedness and almost-gradedness of the algebra under the commutator.  
There is a geometric reason behind. Note that the minimal codimension of 
a strict symplectic subspace of a symplectic vector space is two. 
\end{remark}
\section*{Acknowledgements}
%%%%%%%%%%%%%%%%%%%%
This work  was partially supported  in the frame of 
 the 
%\newline
ESF Research Networking Programmes
{\it Harmonic and Complex Analysis and Applications}, HCAA,
and {\it Interactions of Low-Dimensional Topology and Geometry
with Mathematical Physics}, ITGP.

%\section{References}

\end{document}